\documentclass[11pt,reqno]{amsproc}

\title[High Reynolds numbers]{Remarks on high Reynolds numbers hydrodynamics and the inviscid limit}

\author[P. Constantin]{Peter Constantin}
\address{Department of Mathematics, Princeton University, Princeton, NJ 08544}
\email{const@math.princeton.edu}
\author[V. Vicol]{Vlad Vicol}
\address{Department of Mathematics, Princeton University, Princeton, NJ 08544}
\email{vvicol@math.princeton.edu}

\usepackage[margin=1.25in]{geometry}
\usepackage{amsmath, amsthm, amssymb}
\usepackage{times}
\usepackage{graphicx}

\usepackage{color}
\usepackage[colorlinks=true, pdfstartview=FitV, linkcolor=blue, citecolor=blue, urlcolor=blue]{hyperref}

\theoremstyle{plain}

\theoremstyle{definition}

\numberwithin{equation}{section}

\newcommand{\be}{\begin{equation}}
\newcommand{\ee}{\end{equation}}
\newcommand{\pa}{\partial}
\newcommand{\la}{\label}
\newcommand{\fr}{\frac}
\newcommand{\na}{\nabla}
\newcommand{\ba}{\begin{array}{l}}
\newcommand{\ea}{\end{array}}
\newcommand{\Rr}{{\mathbb R}}
\newcommand{\beg}{\begin}

\begin{document}

\begin{abstract}
We prove that any weak space-time $L^2$ vanishing viscosity limit of  a sequence of strong solutions of Navier-Stokes equations in a bounded domain of $\Rr^2$ satisfies the Euler equation if the solutions' local enstrophies are uniformly bounded. We also prove that $t-a.e.$ weak $L^2$ inviscid limits of solutions of 3D Navier-Stokes equations in bounded domains are weak solutions of the Euler equation if they locally satisfy a scaling property of their second order structure function.  The conditions imposed are far away from boundaries, and wild solutions of Euler equations are not a priori excluded in the limit.  
\hfill \today
\end{abstract}

\noindent

\maketitle

\section{Introduction}
Experimentally verified to a large degree, the two-thirds law and the law of finite energy dissipation are cornerstones of turbulence theory (\cite{frisch}).
The law of finite energy dissipation states the non-vanishing of the rate of dissipation of kinetic energy of turbulent fluctuations per unit mass, in the limit of zero viscosity. 
This can be expressed, if Navier-Stokes equations are used, as
\be
\lim_{\nu\to 0}\nu\langle |\na u(x,t)|^2\rangle = \epsilon>0
\la{eps}
\ee
where $\nu$ is the kinematic viscosity, $u$ is the turbulent velocity fluctuation, $\na$ are spatial gradients, and $\langle\dots\rangle$ is a relevant average.

The two-thirds law states that
\be
\langle |u(x+y,t)-u(x,t)|^2\rangle \sim (\epsilon |y|)^{\fr{2}{3}}
\la{2/3}
\ee
for $|y|$ in the inertial range, that is, in a range of scales 
\be
\eta\le |y|\le L
\la{range}
\ee
where $L$ is a the integral scale of turbulence and $\eta$ is the Kolmogorov dissipation scale,
\be
\eta = \nu^{\fr{3}{4}}\epsilon^{-\fr{1}{4}}.
\la{eta}
\ee
The expressions $s_p(y) = \langle |u(x+y,t)-u(x,t)|^p\rangle$ are called  $p$-th order structure functions, and various turbulence theories argue about scaling properties of the type $s_p\sim |y|^{p\zeta_p}$ in the inertial range. Turbulence is parameterized by the Reynolds number
\be
\fr{UL}{\nu}
\la{re}
\ee
where $U$ is a relevant velocity (for instance average r.m.s velocity). 
Often turbulence is generated at boundaries. Thin boundary layers carry significant changes of momentum or heat. Experimentally, in strong turbulence, these boundary layers detach and  heat and momentum are transported to the bulk of the fluid. Much of the dissipation of kinetic energy takes place in the boundary layers.

An asymptotic  description of the vanishing viscosity limit (the high Reynolds number limit, with $U$ and $L$ bounded) was proposed by Prandtl~\cite{Prandtl}. In it, boundary layers of size proportional to $\sqrt{\nu}$ are attached to boundaries. Outside them, the limit is given by the Euler equations. Inside them, a different equation is valid (the Prandtl equation) and there is matching between the two behaviors at the edges of the boundary layer. If such a description is valid, then zero viscosity limits of solutions of the Navier-Stokes equations inside the domain obey the Euler equations.
 
Much effort has been devoted to validate mathematically turbulence theories and the inviscid limit to the Euler equations. One of the most interesting connections between the two subjects has been made by Kato (\cite{kato}). He proved the equivalence of four statements, for short time, in a regime in which the Euler equations are smooth and conserve energy. These are:

1. The strong convergence in $L^{\infty}(0,T; L^2(\Omega))$
\[
\lim_{\nu\to 0}\sup_{t\in[0,T]}\|u^{NS}(t)-u^E(t)\|_{L^2(\Omega)} = 0.
\]

2. The weak convergence in $L^2(\Omega)$ for {\em{all}} fixed times of the velocity of the Navier-Stokes solution $u^{NS}(t)$ to the velocity of the Euler solution, 
$u^{E}(t)$.

3. The vanishing of the energy dissipation rate
\[
\lim_{\nu\to 0}\nu\int_0^T\|\na u^{NS}(t)\|^2_{L^2(\Omega)}dt = 0.
\]

4. The vanishing of the energy dissipation rate in a very thin boundary layer
of width proportional to $\nu$, $\Gamma_{\nu}$:
\[
\lim_{\nu\to 0}\nu\int_0^T\|\na u^{NS}(t)\|^2_{L^2(\Gamma_\nu)}dt = 0.
\]
The result is a stability result of the Euler path $S^E(t)u_0$, conditioned on assumptions on the viscous dissipation at the boundary. There is a large literature concerned with related or similar conditional strong $L^2$ convergence results (a few examples are ~\cite{BardosTiti07,BardosTiti13, ceiv,ConstantinKukavicaVicol15,Kelliher07, Kelliher08,TemamWang97b,Wang01}).
Some strong $L^2$ unconditional convergence results for short time do exist. They are based on assumptions of real analytic data~\cite{SammartinoCaflisch98b}, or the vanishing of the Eulerian~ initial vorticity in a neighborhood of the boundary~\cite{Maekawa14}. Symmetries can also lead to strong inviscid limits~\cite{Anna15,Kelliher09,LopesMazzucatoLopes08,LopesMazzucatoLopesTaylor08,MazzucatoTaylor08}. All these unconditional results are for short time, close to a smooth solution of Euler equation in laminar situations where energy dissipation rates vanish in the limit. The vast majority of the conditional results are also for short time, close to a smooth solution of Euler equation in laminar situations where energy dissipation rates vanish in the limit, and the conditions involve some uniform property of the Navier-Stokes solutions near the boundary such as bounds on derivatives (like the wall shear stress) or at least some uniform equicontinuity (\cite{ceiv}).

What happens in the bulk for turbulent flows in domains with boundaries is a fundamental open problem. Is there a connection between the Euler equations and the inviscid limit when the limiting energy dissipation rate does not vanish?

Infinite time and the zero vsicosity limit do not commute. This is obvious in the case of unforced Navier-Stokes equations in a smooth regime without boundaries, where the infinite time viscous limits are all zero, and the finite time inviscid limits are conservative smooth Euler solutions. This lack of interchangeability of limits is also true in the forced case. Consider, for example a sequence of solutions of two dimensional, spatially periodic solutions of the Navier-Stokes equations with Kolmogorov forcing $f$,
i.e., forces which are eigenfunctions of the Stokes operator $A$:
\[
Af = \lambda f.
\]
Unique, exact solutions of the Navier-Stokes equations are of the form $u(t) = y(t)f$ with the real valued function $y(t)$ given by
\[
y(t) = y_0e^{-\nu\lambda t} + \fr{1}{\nu \lambda}\left(1-e^{-\nu\lambda t}\right).
\]
Exact solutions of the Euler equations are of the same form, with
\[
y(t) = y_0 + t.
\]
For any finite time the Navier-Stokes solutions converge to the Euler solution $S^{NS}(t)u_0\to S^{E}(t)u_0$ and the {solutions} are bounded as $\nu\to 0$, locally in time. By contrast, the infinite time limit at fixed viscosity is $u(t)\to u_f = \fr{1}{\nu\lambda }f$, and this {sequence} obviously diverges as $\nu\to 0$. Also, the initial data are forgotten in the infinite time limit. If the forcing has odd symmetry, the solutions obey Dirichlet boundary conditions as well.

Because of the lack of interchangeability of limits it is  important to distinguish between the short time zero viscosity limit, the arbitrary finite time limit, and the infinite time limit. 

In this paper we prove two results. They are for arbitrary finite time, and the conditions imposed are far away from boundaries. The results are of weak convergence on subsequences, allowing for non-unique, possibly dissipative Euler limit solutions. 

For {2D} flows at high Reynolds numbers we prove that any $L^2(0,T; L^2(\Omega))$ weak limit of a sequence of strong solutions of Navier-Stokes equations satisfies the Euler equations if interior local enstrophy bounds are uniform in viscosity. No assumptions need to be placed on the behavior of the Navier-Stokes solutions near the boundary. This is not a stability result of an Euler path, but rather a reflection interior good behavior of Navier-Stokes solutions uniform in viscosity. The {{limiting}} Euler solutions inherit interior enstrophy bounds, but the energy dissipation rate might be non-vanishing in the limit of zero viscosity. 

For 3D we prove that if $S^{NS}(t)u_0$ converge weakly in $L^2(\Omega)$ for {{almost all time}}  to a function $u_{\infty}(t)$, and if a second order structure function scaling from above is assumed locally uniformly (like in the two-thirds law, but with any positive exponent), then $u_{\infty}$  satisfies Euler equations. This is different than Kato's condition 2 in that no assumptions are placed on $u_{\infty}$, and all time convergence is not required. In fact, the rate of dissipation of energy need not vanish in the limit, no Euler path is singled out, and the Euler solution may be wild.     

We start by establishing the notation and make preliminary comments. Section 2 is devoted to 2D and section 3 to 3D. A brief discussion concludes the paper.

We consider a bounded open domain $\Omega\subset \Rr^d$ $d=2,3$ with smooth boundary. 
We denote by $u:\Omega\times [0,T)\to \Rr^d$  a solution of the Navier-Stokes equation
\be
\pa_t u -\nu\Delta u + u\cdot\na u + \na p = f
\la{nse}
\ee
in $\Omega$ with
\be
\na\cdot u =0,
\la{divu}
\ee
boundary conditions
\be
u_{\left |\right.  \pa \Omega} =0,
\la{bc}
\ee
and initial data
\be
u_{t=0} = u_0.
\la{ic}
\ee
The velocity $u =  S^{NS}(t)(u_0)$ obviously depends on $\nu>0$, space variable $x\in\Omega$, time variable $t\in [0,T)$, with $T$ possibly infinite, body forces $f$, and initial data $u_0$. 

We discuss weak limits in $L^2(0,T; L^2(\Omega))$. The existence of weak 
limits of solutions of the Navier-Stokes {{equations}} is guaranteed by bounds
\be
\int_0^T\int_{\Omega}|u(x,t)|^2dxdt \le E
\la{enunif}
\ee
which are uniform for the ensemble of solutions. Conversely, if a weak limit
exists for a sequence of functions, (\ref{enunif}) is necessary.  If a sequence $u_n$ converges weakly in $L^2(0,T; L^2(\Omega))$ to a function $u$, it does not follow that $u_n(t)$ converges weakly to $u(t)$ in $L^2(\Omega)$ for almost all $t$, not even on a subsequence. A subsequence of a weakly convergent sequence converges weakly to the same limit, and the subsequence might have some additional properties. In this paper we use this fact to deduce additional information about the weak limits in two dimensions.

We say that function $u\in L^2(0,T;L^2(\Omega))$ is a weak solution of the Euler equations if it is divergence-free and satisfies the Euler equations in the sense of distributions:
\be
(u, \Phi_t) + ((u\otimes u):\na\Phi) + (f,\Phi) = 0
\la{weuler}
\ee
for any $\Phi\in C_0^{\infty}((0,T)\times \Omega)$ which is divergence-free. The notation $M:N$ refers to the trace of the product of the two matrices. This is the distributional form of the incompressible Euler equations 
\[
\pa_t u + u\cdot\na u + \na p = f, \quad \na\cdot u = 0
\]
forced by $f$. No boundary conditions nor initial data are part of the distributional formulation.

In order to verify that the limit $u$ of a weakly convergent sequence $u_n \in L^2(0,T; L^2(\Omega))$ of solutions of Navier-Stokes equations  satisfies the Euler equations, it is enough to prove the convergence
\be
N_{\Phi}(u_n)\to N_{\Phi}(u)
\la{nphilim}
\ee
for any fixed, divergence-free test function $\Phi\in C_0^{\infty}((0,T)\times \Omega)$ where  
\be
N_\Phi(u) = \int_0^T\int_{\Omega}(u\otimes u):\na \Phi\, dx dt.
\la{Nphins}
\ee
This is true of course only if we assume that the forces driving the Navier-Stokes equations converge weakly in $L^2$ to $f$. Then the linear terms (viscous term, time derivative term) and the forcing terms obviously converge. 

In two dimensions, we use a vorticity formulation of the equations.

\section{2D}

We consider the vorticity $\omega = \pa_1u_2-\pa_2u_1 = \na^{\perp}\cdot u$. For solutions of the 2D Navier-Stokes equations the vorticity obeys
\be
\pa_t \omega + u\cdot\na \omega -\nu\Delta \omega = g = \na^{\perp}\cdot f.
\la{omegaeq}
\ee
We recall that the velocity is obtained from a stream function 
\be
u = \nabla^{\perp}\psi
\la{psi}
\ee
and that 
\be
\omega = \Delta\psi
\la{omegapsi}
\ee
and therefore 
\be
\Delta u = \na^{\perp}\omega
\la{deltauomega}
\ee
holds. Note that $-\Delta u$ is not the Stokes operator applied to $u$.
The identity
\be
\int_{\Omega}|\omega(x,t)|^2dx = \int_{\Omega}|\na u(x,t)|^2dx
\la{enstgrad}
\ee
is true in view of the boundary conditions (\ref{bc}). Indeed:
\[
\ba
\int_{\Omega} \pa_j u_i(x,t)) \pa_j u_i(x,t)dx  = -\int_{\Omega} u(x,t)\cdot\Delta u(x,t)dx \\
=- \int_{\Omega}u(x,t)\cdot\na^{\perp}\omega(x,t)dx = \int_{\Omega}|\omega(x,t)|^2dx.
\ea
\]
In the integrations by parts we used only (\ref{bc}). For (\ref{deltauomega}),
the fact that $u(t)\in H^2(\Omega)$ for almost all time is true  because $u$ is a strong solution of NSE but in fact the equality (\ref{enstgrad}) 
is true for all divergence-free $u\in H_0^1(\Omega)$, by approximation.

\beg{thm} Let $\Omega\subset\Rr^2$ be a bounded open set with smooth boundary. Let $u_n$ be a sequence of solutions of Navier-Stokes equations with viscosities $\nu_n\to 0$. We assume that the solutions are driven by forces $f_n\in H^1(\Omega)$ that are uniformly bounded in $H^1(\Omega)$ and converge weakly in $H^1(\Omega)$ to $f$. We assume that the initial data $u_n(0)$ are divergence-free and are uniformly bounded in $H^1_0(\Omega)$.  Let $K$ be a compact, $K\subset\subset \Omega$. We assume that there exists a constant ${\mathcal{E}}_K$
 which might depend on $H^1(\Omega)$ norms of initial data and $f$, on $K$ and $T$, but is  {\em{independent of viscosity}}, such that
\be
\sup_{0\le t\le T}\int_K|\omega_n(x,t)|^2dx\le {\mathcal E}_K
\la{locenst}
\ee
where $\omega_n = \nabla^{\perp}\cdot u_n$ are the vorticities.

Then any weak limit in $L^2(0,T; L^2(\Omega))$ of the sequence $u_n$, $u_{\infty}$, is a weak solution of the Euler equations
\be
\pa_t \omega_{\infty} + u_{\infty}\cdot\na\omega_{\infty} = g =\na^{\perp}\cdot f
\la{vortinftyeq}
\ee
with $\omega_{\infty} = \na^{\perp}\cdot u_{\infty}$. The solution has bounded energy,
\be
u_{\infty}\in L^{\infty}(0,T; L^2(\Omega)).
\la{bene}
\ee
Moreover, for any compact $K\subset\subset \Omega$ there exists a constant $\widetilde{\mathcal E}_K$ such that
\be
\sup_{t\in[0,T]}\int_K |\omega_{\infty}(x,t)|^2dx \le \widetilde{\mathcal E}_K
\la{loco}
\ee
holds. 
\end{thm}
\noindent{\bf Proof.} The fact that $u_{\infty}$ has bounded energy is a simple consequence of the fact that under our conditions the sequence $u_n$ is bounded in time in $L^2(\Omega)$. Indeed, for time intervals $I$, $\chi_I(t) u_n$ converge weakly in $L^2(0,T; L^2(\Omega))$ to $\chi_I(t)u_{\infty}$ where $\chi_I$ is the indicator function of $I$. Thus
\[
\int_I\|u_{\infty}(t)\|^2_{L^2(\Omega)}dt \le \lim\inf_{n\to\infty}\int_I\|u_n(t)\|^2_{L^2(\Omega)}dt \le C|I|,
\]
and (\ref{bene}) follows.
In order to prove that $u_{\infty}$ solves (\ref{vortinftyeq}), we consider the nonlinear term, which is the only term whose behavior is in question. 
We take a compactly supported test function $\Phi\in C_0^{\infty}((0,T)\times\Omega)$ whose  support is a compact $L\subset [t_1, T_1]\times K_1$ with $K_1\subset\subset \Omega$ compact, and $0<t_1<T_1<T.$ We consider a larger compact $K\subset\subset\Omega$ such that $K_1$ is included in the interior of $K$ and a slightly larger time interval $[t_0, T_0]$ with $0<t_0<t_1$ and $T_1<T_0<T$. Let us also take a function $\chi_0\in C_0^{\infty}((t_0, T_0)\times K)$ which is identically 1 on a neighborhood of $[t_1, T_1]\times K_1$. We consider the sequence $\chi_0 u_n$. Because $u_n$ are uniformly bounded in $L^2(\Omega)$ it is clear  in view of (\ref{locenst}) that $\na^{\perp}(\chi_0 u_n)$ is a bounded sequence in $L^{\infty}(0,T; L^2(\Omega))$. Because $\na\cdot (\chi_0 u_n)$ is bounded in $L^{\infty}(0,T, L^2(\Omega))$ as well, it follows that $\chi_0 u_n$ is bounded in $L^{\infty}(0,T; H_0^1(\Omega))$. In order to use a Aubin-Lions lemma, and obtain some uniform control on time derivatives it is best to take the curl of the equation, because the vorticity equation is local. 
The equations obeyed by the vorticities are
\be
\pa_t \omega_n + u_n\cdot\na\omega_n -\nu\Delta \omega_n = g_n.
\la{vorteqn}
\ee
We consider now another cutoff function $\chi$ which is still equal identically to 1 on a neighborhood of $[t_1,T_1]\times K_1$ but whose compact support is included in the region where $\chi_0$ is identically 1.
We multiply by $\chi$ and consider the sequence $w_n = \chi\omega_n$. In view of  (\ref{locenst}) $w_n$ is bounded in $L^{\infty}(0,T; L^2(\Omega))$. We use the equation (\ref{vorteqn}) to examine $\pa_t w_n$. 

The sequence $\chi \nu\Delta \omega_n$ is bounded in $L^{\infty}(0,T; H^{-2}(\Omega))$, where $H^{-2}(\Omega)$ is the dual of $H_0^2(\Omega)$ because of (\ref{locenst}). The terms $\pa_t \chi \omega_n$ and $\chi g_n$ are bounded in a better space, $L^{\infty}(0,T; L^2(\Omega))$. It is well known that the term $u_n\cdot\na\omega_n$ is a second derivative. Indeed, dropping the subscript $n$ for a moment in order to avoid confusion, 
\be 
u\cdot\na \omega = \pa_1\pa_2(u_2^2-u_1^2) + (\pa_1^2-\pa_2^2)(u_1u_2)
\la{unaotwo}
\ee
where now $u_1, u_2$ are components of the vector $u_n$. Therefore, because on the support of $\chi$ we have that $u_n = \chi_0 u_n$, it follows that the term
$\chi u_n\cdot\na\omega_n$ is bounded in $L^{\infty}(0,T; H^{-2}(\Omega))$. Indeed, using the continuous embedding $H_0^1(\Omega)\subset L^4(\Omega)$ we have that $\chi_0 u_n$ are uniformly bounded in $L^4(\Omega)$, and thus, after peeling off the two derivatives of $u_n\cdot\na \omega_n$ we are left with functions that are bounded uniformly in $L^{\infty}(0,T; L^2(\Omega))$. By the Aubin-Lions lemma with spaces $L^2(\Omega)\subset\subset H^{-1}(\Omega)\subset H^{-2}(\Omega)$ (\cite{Lions}) we have that the sequence $w_n$ has a strongly convergent subsequence (relabelled $w_n$) in $L^2(0,T; H^{-1}(\Omega))$. 
More precisely, we have that $\Lambda_D^{-1}w_n$ converges strongly in $L^2(0,T; L^2(\Omega))$ to a function $v$, where $\Lambda_D = (-\Delta)^{\fr{1}{2}}$ with $-\Delta$ the Laplacian with homogeneous Dirichlet boundary conditions in $\Omega$. It is well known that $\Lambda_D:H_0^1(\Omega)\to L^2(\Omega)$ is an isometry. 
Taking a test function $\Psi(x,t)$ we have that
\[
\int_0^T\int_{\Omega}\Lambda_D^{-1}w_n(x,t)\Psi(x,t) = -\int_0^T\int_{\Omega}
u_n(x,t)\cdot \na^{\perp} {(\chi \Lambda^{-1}_D\Psi)}(x,t)dxdt.
\]
We pass to the limit in both sides, noting that $\na^{\perp}{(\chi \Lambda^{-1}_D\Psi)}(x,t)$ is an allowed test function because and it belongs to $L^2(0,T; L^2(\Omega))$, in view of the boundedness of the Riesz transforms $R_D = \nabla\Lambda_D^{-1}$ in $L^2(\Omega)$.
It follows that
\[
v= \Lambda_D^{-1}(\chi\na^{\perp}\cdot u_{\infty}) = \Lambda_D^{-1}\chi\omega_{\infty}
\]
Moreover, because $\|\Lambda_D^{-1}w_n(t)\|_{L^2(\Omega)}^2$ converge strongly in $L^1(0,T)$ there is a subsequence, relabelled by $n$ such that $\Lambda_D^{-1}w_n(t)$ converges strongly in $L^2(\Omega)$ for almost all $t\in [0,T]$ to $v(t)$. Testing with a test function $\phi$ we have that, on one hand 
\[
\left |\int_{\Omega}\phi w_n(t)dx\right| =\left |\int_{\Omega}\chi\omega_n(t)\phi dx\right|\le {\mathcal E}_K^{\fr{1}{2}}\|\phi\|_{L^2(\Omega)},
\]
and on the other
\[
\int_{\Omega}(\Lambda_D \phi)v(t)dx = \lim_{n\to\infty}\int_{\Omega}\Lambda_D\phi (\Lambda_D^{-1}w_n(t))dx
\]
holds for almost all $t$. Thus
\[
\int_{\Omega}\left| \chi\na^{\perp}\cdot u_{\infty}(x,t)\right |^2dx \le{\mathcal E}_K
\]
holds a.e. in time, that is
\be
\|\chi \omega_{\infty}\|_{L^{\infty}(0,T; L^2(\Omega))}\le {\mathcal E}_K^{\fr{1}{2}}.
\la{omeglocb}
\ee
We now pass to another subsequence, relabelled again by $n$ such that $\chi_0 u_n$ converges weakly in $L^2(0,T; H_0^1(\Omega))$. It is obvious that the limit is $\chi_0 u_{\infty}$.
We consider the nonlinear term
\[
\int_0^T\int_\Omega (u_n\cdot\na \Phi)\omega_ndxdt = N_{\Phi}(n)
\]
Because on the support of $\Phi$ we have $u_n = \chi_0 u_n$ and $\omega_n = w_n$, this integral is the duality pairing between a weakly convergent sequence in
$L^2(0,T; H_0^1(\Omega))$, namely $\chi_0u_n\cdot\na\Phi$ and a strongly convergent sequence in $L^2(0,T; H^{-1}(\Omega))$, namely $w_n$:
\[
N_{\Phi}(n) = \int_0^T\int_{\Omega} \Lambda_D(\chi_0u_n\cdot\na\Phi)\Lambda_D^{-1}w_ndx dt.
\]
Therefore $N_\Phi(n)$ convergences as the scalar product between weakly convergent and strongly convergent $L^2$ functions, and  
\be
\lim_{n\to\infty} N_{\Phi}(n) = \int_0^T\int_{\Omega}(\chi_0u_{\infty}\cdot\na\Phi)\chi\omega_{\infty}dxdt = \int_0^T\int_{\Omega}(u_{\infty}\cdot\na \Phi)\omega_{\infty}dxdt
\la{Nphi}
\ee
\beg{rem} We note that from the proof it follows that $\chi \omega_n(t)$ converge weakly in $L^2$ to $\chi\omega_{\infty}(t)$ on a subsequence, for almost all $t$. No convergence is implied for $u_n(t)$ in $L^2(\Omega)$: the global behavior may depend on viscosity. 
\end{rem}

\section{3D}
We consider families of solutions of Navier Stokes equations in a bounded domain $\Omega\subset\Rr^3$.
\beg{thm} Let $u_n$ be a sequence of weak solutions of the Navier-Stokes equations
\be
\pa_t u_n + u_n\cdot\na u_n -\nu_n\Delta u_n + \na p_n = f_n
\la{nsn}
\ee
with $\na\cdot u_n =0$, $f_n$ bounded in $L^2(0,T; L^2(\Omega))$, converging weakly to $f$, $u_n(0)$ divergence-free and bounded in $L^2(\Omega)$ and $\nu_n\to 0$. We assume that for any $K\subset\subset \Omega$ there exists a constant $E_K$, a constant $\epsilon>0$ and a constant $\zeta_2>0$ such that
\be
\sup_{n}\int_0^T\int_K |u_n(x+y,t)-u_n(x,t)|^2dxdt \le E_K |y|^{2\zeta_2}
\la{s2}
\ee
holds for $|y|< dist(K,\pa\Omega)$ in the inertial range
\be
|y|\ge \epsilon^{-\fr{1}{4}}\nu_n^{\fr{3}{4}} = \eta(n)
\la{ir}
\ee
Assume that $u_n(t)$ converge weakly in $L^2(\Omega)$ to $u_{\infty}(t)$ for almost all $t\in (0,T)$. Then $u_{\infty}$ is a weak solution of the Euler equations.
\end{thm}
\beg{rem} The domain need not be bounded. Local uniform energy bounds are enough. In exterior domains local uniform (in viscosity) energy bounds can be obtained for suitable weak solutions. 
\end{rem}
\beg{rem} Obviously, the scaling assumption (\ref{s2}) does not imply regularity, because it is limited to $y$ bounded away from zero. Also, the exact Kolmogorov form of $\eta(n)$ is not needed. All that is used is that $\eta(n)$ converges to zero as $n\to\infty$. Finally, the power law behavior is not needed either, any uniform modulus of continuity can be used instead.
\end{rem}
\noindent{\bf Proof.} We consider a nonnegative smooth function $j(z)$ supported in the annulus $1<|z|<2$ and with integral equal to 1, $\int_{\Rr^3} j(z)dz =1$. We assume also that $j(-z)=j(z)$. We fix a compact $K\subset\subset \Omega$ and denote, for a function $u$, for $x\in K$, and $2r<dist (K, \pa\Omega)$,  
\be
u_r(x) = \int_{1\le |z|\le 2}u(x-rz)j(z)dz.
\la{ur}
\ee
We note the identity  (see \cite{cet})
\be
(uv)_r(x)- u_r(x)v_r(x) = \rho_r(u,v)(x)
\la{uvrho}
\ee
with
\be
\rho_r(u,v)(x) = \int_{1\le |z|\le 2} j(z)(u(x-rz)-u(x))(v(x-rz)-v(x))dz  - (u(x)-u_r(x))(v(x)-v_r(x)).
\la{rhor}
\ee
Let us take a divergence-free test function $\Phi(x,t)\in (C_0^{\infty}((0,T)\times\Omega))^3$ and investigate the behavior of
the nonlinear term 
\be
N_\Phi(n) = \int_0^T\int_{\Omega}(u_n\otimes u_n):\na \Phi\, dx dt.
\la{Nphin}
\ee
We take a compact $K$ such that the support of $\Phi$ is included in 
$[t_0, T_0]\times K$, where $0<t_0<T_0<T$.
Let us start by noting that 
\[
\left| \int_0^T\int_{\Omega}(u_n\otimes u_n) : (\na\Phi -(\na\Phi)_r)\,dxdt\right| \le C_{\Phi}r\|u_n\|^2_{L^2(0,T; L^2(K))}
\]
and thus
\be
\left| \int_0^T\int_{\Omega}(u_n\otimes u_n) : (\na\Phi -(\na\Phi)_r)\,dxdt\right| \le rC_{\Phi}E
\la{nr}
\ee
with $E$ a uniform bound on the local time average of energy,
\be
\|u_n\|^2_{L^2(0,T; L^2(K))}\le E.
\la{ene}
\ee
Then, we note, using $j(z)= j(-z)$ that
\be
\int_0^T\int_{\Omega}u_n\otimes u_n : (\na\Phi)_r \,dxdt = 
\int_0^T\int_{\Omega}(u_n\otimes u_n)_r : \na\Phi \,dxdt.
\la{switchr}
\ee 
Now we use the identity (\ref{uvrho})
\be
\ba
\int_0^T\int_{\Omega}(u_n\otimes u_n)_r: (\na\Phi) \,dxdt =\\
\int_0^T\int_{\Omega}(u_n)_r\otimes (u_n)_r : \na\Phi \,dxdt +
\int_0^T\int_{\Omega}\rho_r(u_n, u_n) : (\na\Phi) \,dxdt.
\ea
\la{sumr}
\ee
We take $n$ large enough so that $\eta(n)\le r$.
We estimate the second term using the assumption (\ref{s2}):
\[
\left|\int_0^T\int_{\Omega}\rho_r(u_n, u_n) : (\na\Phi) \,dxdt\right|\le
C_{\Phi}E_K r^{2\zeta_2}
\]
Therefore we have proved so far that
\be
\left| N_{\Phi}(n)- \int_0^T\int_{\Omega}(u_n)_r\otimes (u_n)_r : \na\Phi \,dxdt \right |\le C_\Phi(E_K r^{2\zeta_2} + E r)
\la{inter}
\ee
holds for $n$ large enough, {depending on $r$}. 

We note that if $u_n(t)$ converges weakly in $L^2(\Omega)$ to $u_{\infty}(t)$, then $(u_n(t))_r(x)$
converges pointwise in $K$ to $(u_{\infty}(t))_r(x)$ at fixed $r$, just because it is the scalar product
\[
(u_n(t))_r(x) = \int_{x-rA}u_n(y,t)j\left(\fr{x-y}{r}\right)r^{-3}dy
\]
where we denoted by $A$ the annulus $1\le |z|\le 2$. For $x\in K$ we observe that $x-rA\subset\subset \Omega$. We also have
\[
|(u_n(t))_r(x)|\le Cr^{-\fr{3}{2}}\|u_n(t)\|_{L^2(\Omega)}\le Er^{-\fr{3}{2}}
\]
which allows us to use the dominated convergence theorem and pass to the limit.
{By the triangle inequality we obtain}
\be
\left| N_{\Phi}(n)- \int_0^T\int_{\Omega}(u_{\infty})_r\otimes (u_{\infty})_r : \na\Phi \,dxdt \right |\le C_\Phi(E_K r^{2\zeta_2} + E r)
\la{interm}
\ee
for $n$ sufficiently large, {depending on $r$}, and $r$ small enough. Now we use the identity
(\ref{uvrho}) in reverse, and the fact that translation is strongly continuous in $L^2$ to deduce that $(u_{\infty})_r$ converges to $u_\infty$ strongly in $L^2(0,T; L^2(K))$. Thus, given $\delta>0$ we can choose $r=r(\delta)$ small enough so that
\[
\int_0^T\int_K|u_{\infty}-(u_{\infty})_r|^2dxdt \le \fr{\delta}{2}
\]
and, using (\ref{interm}) and making sure that $r$ is small enough that
\[
 C_\Phi(E_K r^{2\zeta_2} + E r)\le \fr{\delta}{2}
\]
holds as well, we obtain
\be
\left| N_{\Phi}(n)- \int_0^T\int_{\Omega}(u_{\infty})\otimes (u_{\infty}) : \na\Phi \,dxdt \right |\le \delta
\la{interme}
\ee
for $n$ large enough. We have thus
\be
\lim_{n\to\infty} N_{\Phi}(n) = \int_0^T\int_{\Omega}(u_{\infty})\otimes (u_{\infty}) : \na \Phi \,dxdt,
\la{final}
\ee 
and this concludes the proof.
\beg{rem} It is possible to remove the assumption of almost all time $L^2(\Omega)$ convergence, and replace it with the weak convergence in $L^2(0,T; L^2(\Omega))$, at the price of demanding space-time second order structure function scaling. 
\be
\int_0^T\int_K |u_n(x+y,t+s)-u_n(x,t)|^2dxdt \le E_K(|y|^{2\zeta_2} + |s|^\beta)
\la{spacetimesc}
\ee
for $\eta(n)\le |y|< dist(K; \pa\Omega)$, $t+s\in [0,T]$, $|s|\ge \tau(n)$, 
$\tau(n)\to 0$, and  $\beta>0$. 
The proof is the same, we translate in space-time. If $\tau(n)=0$ the requirement is strong, it implies the sequence bounded in $C^{\beta}(0,T; L^2(\Omega))$, and in particular the $L^2(\Omega)$ convergence on each time slice.
\end{rem}
\beg{rem} By Fatou's lemma in time and our assumptions, it follows that the limit solution of Euler equations satisfies the local bounds
\be
\int_0^T\int_K |u_{\infty}(x+y,t)-u_{\infty}(x,t)|^2dxdt \le E_K |y|^{2\zeta_2}
\la{s2fty}
\ee
for $|y|<dist (K, \pa \Omega)$ and any compact $K\subset\subset\Omega$.
\end{rem}
\section{Discussion}
The vanishing of the dissipation rate follows from weak convergence in $L^2(\Omega)$ for all times only if the Euler equation is conservative. We proved results of emergence of weak, possibly dissipative solutions of Euler equations in 3D if the ensemble of Navier-Stokes solutions obeys a local-in-space but uniform in the ensemble second order structure function scaling from above. In two dimensions, we proved the emergence of weak solutions form arbitrary families of strong solutions of Navier-Stokes equations with uniform interior (local) enstrophy bounds. There might be dissipative solutions among them, although an example is not available at this time.

\vspace{.2in}
\noindent{\bf{Acknowledgment.}} The research of PC is partially funded by NSF grant DMS-1209394, and the research if VV is partially funded by NSF grant DMS-1652134.

\end{document}